\documentclass[12pt,a4paper]{article}
\usepackage{amssymb}

\usepackage{graphicx,amssymb,amsfonts,epsfig,amsthm,a4,amsmath,url}
\usepackage[latin1]{inputenc}

\newtheorem{ThmIntro}{Theorem}
\newtheorem{CorIntro}[ThmIntro]{Corollary}
\newtheorem{PropIntro}[ThmIntro]{Proposition}

\newtheorem{thm}{Theorem}[section]
\newtheorem{cor}[thm]{Corollary}
\newtheorem{lem}[thm]{Lemma}

\newtheorem{prop}[thm]{Proposition}
\theoremstyle{definition}
\newtheorem{defn}[thm]{Definition}

\newtheorem{que}[thm]{Question}

\theoremstyle{remark}
\newtheorem{rem}[thm]{Remark}
\newtheorem{ex}[thm]{Example}

\numberwithin{equation}{section}

\newcommand{\Z}{\mathbf{Z}}

\newcommand{\N}{\mathbf{N}}
\newcommand{\R}{\mathbf{R}}

\newcommand{\Q}{\mathbf{Q}}

\newcommand{\HH}{\mathcal{H}}

\newcommand{\bpr}{\noindent \textbf{Proof}: ~}

\newcommand{\epr}{~$\blacksquare$}

\newcommand{\eps}{\varepsilon}

\newcommand{\Var}{\textnormal{Var}}
\newcommand{\Supp}{\textnormal{Supp}}
\newcommand{\Exp}{\textnormal{Exp}}

\title{Asymptotic isoperimetry on groups and uniform embeddings into Banach spaces.}
\author{Romain Tessera}
\date{\today}

\begin{document}

\baselineskip=16pt

\maketitle

\begin{abstract}
We characterize the possible asymptotic behaviors of the
compression associated to a uniform embedding into some
$L^p$-space, with $1<p<\infty$, for a large class of groups
including connected Lie groups with exponential growth and
word-hyperbolic finitely generated groups. In particular, the
Hilbert compression rate of these groups is equal to $1$. This
also provides new and optimal estimates for the compression of a
uniform embedding of the infinite 3-regular tree into some
$L^p$-space. The main part of the paper is devoted to the explicit
construction of affine isometric actions of amenable connected Lie
groups on $L^p$-spaces whose compressions are asymptotically
optimal. These constructions are based on an asymptotic lower
bound of the $L^p$-isoperimetric profile inside balls. We compute
the asymptotic of this profile for all amenable connected Lie
groups and for all $1\leq p<\infty$, providing new geometric
invariants of these groups. We also relate the Hilbert compression
rate with other asymptotic quantities such as volume growth and
probability of return of random walks.
\end{abstract}

\tableofcontents



\section{Introduction}

The study of uniform embeddings of locally compact groups into
Banach spaces and especially of those associated to proper affine
isometric actions plays a crucial role in various fields of
mathematics ranging from K-theory to geometric group theory.
Recall that a locally compact group is called a-T-menable if it
admits a proper affine action by isometries on a Hilbert space
(for short: a proper isometric Hilbert action). An amenable
$\sigma$-compact locally compact group is always a-T-menable
\cite{CCJJV}; but the converse is false since for instance
non-amenable free groups are a-T-menable. However, if a locally
compact, compactly generated group $G$ admits a proper isometric
Hilbert action whose compression $\rho$ satisfies
$$\rho(t)\succ t^{1/2},$$ then $G$ is amenable\footnote{This was
proved for finitely generated groups in \cite{GK}. In
\cite{coteva}, we give a shorter argument that applies to all
locally compact compactly generated groups. }. On the other hand,
in \cite{coteva}, we prove that non-virtually abelian polycyclic
groups cannot have proper isometric Hilbert actions with linear
compression. These results motivate a systematic study of the
possible asymptotic behaviors of compression functions, especially
for amenable groups.

In this paper, we ``characterize" the asymptotic behavior of the
$L^p$-compression, with $1<p<\infty$, for a large class of groups
including all connected Lie groups with exponential growth. Some
partial results in this direction for $p=2$ had been obtained in
\cite{GK} and \cite{BS} by completely different methods.


\subsection{$L^p$-compression: optimal estimates}

Let us recall some basic definitions. Let $G$ be some locally
compact compactly generated group. Equip $G$ with the word length
function $|\cdot|_S$ associated to a compact symmetric generating
subset $S$ and consider a uniform embedding $F$ of $G$ into some
Banach space. The compression $\rho$ of $F$ is the nondecreasing
function defined by
$$\rho(t)=\inf_{|g^{-1}h|_S\geq t}\|
F(g)-F(h)\|.$$

Let $f,g: \R_+\to\overline{\R}_+$ be nondecreasing, nonzero
functions. We write respectively $f\preceq g$, $f\prec g$ if there
exists $C>0$ such that $f(t)=O(g(Ct))$, resp. for all $c>0$,
$f(t)=o(g(ct))$ when $t\to \infty$. We write $f\approx g$ if both
$f\preceq g$ and $g\preceq f$. The asymptotic behavior of $f$ is
its class modulo the equivalence relation $\approx$.

Note that the asymptotic behavior of the compression of a uniform
embedding does not depend on the choice of $S$.

In the sequel, a $L^p$-space denotes a Banach space of the form
$L^p(X,m)$ where $(X,m)$ is a measure space. A
$L^p$-representation of $G$ is a continuous linear $G$-action on
some $L^p$-space. Let $\pi$ be a isometric $L^p$-representation of
$G$ and consider a $1$-cocycle $b\in Z^1(G,\pi)$, or equivalently
an affine isometric action of $G$ with linear part $\pi$: see the
preliminaries for more details. The compression of $b$ is defined
by
$$\rho(t)=\inf_{|g|_S\geq t}\| b(g)\|_p.$$

In this paper, we mainly focus our attention on groups in the two
following classes.

\noindent Denote $(\mathcal{L})$ the class of groups including
\begin{enumerate}
\item polycyclic groups and connected amenable Lie groups;

\item semidirect products
$\Z[\frac{1}{mn}]\rtimes_{\frac{m}{n}}\Z$, with $m,n$ co-prime
integers with\footnote{This condition garanties that the group is
compactly generated.} $|mn|\ge 2$ (if $n=1$ this is the
Baumslag-Solitar group $BS(1,m)$); semidirect products
$\left(\R\oplus\bigoplus_{p\in
P}\Q_{p}\right)\rtimes_{\frac{m}{n}}\Z$ with $m,n$ coprime
integers and $P$ a finite set of primes (possibly infinite,
$\Q_{\infty}=\R$) dividing $mn$;

\item wreath products $F\wr\Z$ for $F$ a finite group.
\end{enumerate}

\noindent Denote $(\mathcal{L'})$ the class of groups including
groups in the class $(\mathcal{L})$ and
\begin{enumerate}
\item connected Lie groups and their cocompact lattices;

\item irreducible lattices in semisimple groups of rank $\geq 2$;

\item hyperbolic finitely generated groups.
\end{enumerate}

Let $\mu$ be a left Haar measure on the locally compact group $G$
and write $L^p(G)=L^p(G,\mu)$. The group $G$ acts by isometry on
$L^p(G)$ via the left regular representation $\lambda_{G,p}$
defined by
$$\lambda_{G,p}(g)\varphi=\varphi(g^{-1}\cdot).$$

\begin{ThmIntro}\label{thm1}
Fix some $1\leq p<\infty$. Let $G$ be a group of the class
$(\mathcal{L})$ and let $f$ be an increasing function $f:\R_+\to
\R_+$ satisfying
\begin{equation}\tag{$C_p$}\int_1^{\infty}\left(\frac{f(t)}{t}\right)^p\frac{dt}{t}<\infty.
\end{equation}
Then there exists a $1$-cocycle $b\in Z^1(G,\lambda_{G,p})$ whose
compression $\rho$ satisfies
$$\rho\succeq f.$$
\end{ThmIntro}


\begin{CorIntro}\label{thm2}
Fix some $1\leq p<\infty$. Let $G$ be a group of the class
$(\mathcal{L'})$ and let $f$ be an increasing function $f:\R_+\to
\R_+$ satisfying Property $(C_q)$, with $q=\max\{p,2\}$. Then
there exists a uniform embedding of $G$ into some $L^p$-space
whose compression $\rho$ satisfies
$$\rho\succeq f.$$
\end{CorIntro}

Let us sketch the proof of the corollary. First, recall
\cite[III.A.6]{W} that for $1\leq p\leq 2$, $L^2([0,1])$ is
isomorphic to a subspace of $L^p([0,1])$. It is thus enough to
prove the theorem for $2\leq p<\infty.$ This is an easy
consequence of Theorem~\ref{thm1} since every group of class
$(\mathcal{L'})$ quasi-isometrically embeds into a group of
$(\mathcal{L})$. Indeed, any connected Lie group admits a closed
cocompact connected solvable subgroup. On the other hand,
irreducible lattices in semisimple groups of rank $\geq 2$ are
quasi-isometrically embedded \cite{LMR}. Finally, any hyperbolic
finitely generated group quasi-isometrically embeds into the real
hyperbolic space $\mathbf{H}^n$ for $n$ large enough \cite{BoS}
which is itself quasi-isometric to $\text{SO}(n,1)$.

The particular case of nonabelian free groups, which are
quasi-isometric to $3$-regular trees, can also be treated by a
more direct method. More generally that method applies to any
simplicial\footnote{By simplicial, we mean that every edge has
length $1$.} tree with possibly infinite degree.
\begin{ThmIntro}\textnormal{(see Theorem~\ref{sphericalprop})}
Let $T$ be a simplicial tree. For every increasing function
$f:\R_+\to\R_+$ satisfying
\begin{equation}\tag{$C_p$}\quad
\int_1^{\infty}\left(\frac{f(t)}{t}\right)^p\frac{dt}{t}<\infty,
\end{equation}
there exists a uniform embedding $F$ of $T$ into $\ell^p(T)$ with
compression $\rho\succeq f.$
\end{ThmIntro}
\begin{rem} In \cite{BuSc,BuSc'}, it is shown that real hyperbolic spaces
and word hyperbolic groups quasi-isometrically embed into finite
products of (simplicial) trees. Thus the restriction of
Corollary~\ref{thm2} to word hyperbolic groups and to simple Lie
groups of rank 1 can be deduced from
Proposition~\ref{sphericalprop}. Nevertheless, not every connected
Lie group quasi-isometrically embeds into a finite product of
trees. Namely, a finite product of trees is a CAT(0) space, and in
\cite{Pauls} it is proved that a non-abelian simply connected
nilpotent Lie group cannot quasi-isometrically embed into any
CAT(0) space.
\end{rem}

\begin{ThmIntro}\label{finitetreethm}
Let $T_N$ be the binary rooted tree of depth $N$. Let $\rho$ be
the compression of some $1$-Lipschitz map from $T_N$ to some
$L^p$-space for $1< p<\infty$. Then there exists $C<\infty$,
depending only on $p$, such that
\begin{eqnarray*}
\int_{1}^{2N}\left(\frac{\rho(t)}{t}\right)^q\frac{dt}{t}\leq C,
\end{eqnarray*}
where $q=\max\{p,2\}$.
\end{ThmIntro}
 This result is a strengthening of
\cite[Theorem 1]{Bourgain}; see also Corollary~\ref{Bourgainthm}.
As a consequence, we have
\begin{CorIntro}
Assume that the 3-regular tree quasi-isometrically embeds into
some metric space $X$. Then, the compression $\rho$ of any uniform
embedding of $X$ into any $L^p$-space for $1<p<\infty$ satisfies
$(C_q)$ for $q=\max\{p,2\}.$
\end{CorIntro}
In \cite[Theorem 1.5]{BenSc}, it is proved that the 3-regular tree
quasi-isometrically embeds into any graph with bounded degree and
positive Cheeger constant (e.g. any non-amenable finitely
generated group). On the other hand, in a work in preparation with
Cornulier \cite{QI}, we prove that finitely generated linear
groups with exponential growth, and finitely generated solvable
groups with exponential growth admit quasi-isometrically embedded
free non-abelian sub-semigroups. Together with the above
corollary, they lead to the optimality of Theorem~\ref{thm1}
(resp. Corollary~\ref{thm2}) when the group has exponential growth
and when $2\leq p<\infty$ (resp. $1<p<\infty$).

\begin{CorIntro}
Let $G$ be a finitely generated group with exponential growth
which is either virtually solvable or non-amenable. Let $\varphi$
be a uniform embedding of $G$ into some $L^p$-space for
$1<p<\infty$. Then its compression $\rho$ satisfies Condition
$(C_q)$ for $q=\max\{p,2\}.$
\end{CorIntro}

\begin{CorIntro}
Let $G$ be a group of class $(\mathcal{L'})$ with exponential
growth. Consider an increasing map $f$ and some $1< p<\infty$;
then $f$ satisfies Condition $(C_q)$ with $q=\max\{p,2\}$ if and
only if there exists a uniform embedding of $G$ into some
$L^p$-space whose compression $\rho$ satisfies $\rho\succeq f$.
\end{CorIntro}
Note that the 3-regular tree cannot uniformly embed into a group
with subexponential growth. So the question of the optimality of
Theorem~\ref{thm1} for non-abelian nilpotent connected Lie groups
remains open.

\

\noindent{\bf About Condition $(C_p)$.} First, note that if $
p\leq q$, then $(C_p)$ implies $(C_q)$: this immediately follows
from the fact that a nondecreasing function $f$ satisfying $(C^p)$
also satisfies $f(t)/t=O(1).$

Let us give examples of functions $f$ satisfying Condition
$(C_p)$. Clearly, if $f$ and $h$ are two increasing functions such
that $f\preceq h$ and $h$ satisfies $(C_p)$, then $f$ satisfies
$(C_p)$. The function $f(t)=t^a$ satisfies $(C_p)$ for every $a<1$
but not for $a=1$. More precisely, the function
$$f(t)=\frac{t}{(\log t)^{1/p}}$$
does not satisfy $(C_p)$ but
$$f(t)=\frac{t}{((\log t)(\log\log t)^a)^{1/p}}$$
satisfies $(C_p)$ for every $a>1$. In comparison, in \cite{BS},
the authors construct a uniform embedding of the free group of
rank $2$ into a Hilbert space with compression larger than
$$\frac{t}{((\log t)(\log\log t)^2)^{1/2}}.$$
As $t/(\log t)^{1/p}$ does not satisfy $(C_p)$, one may wonder if
$(C_p)$ implies
$$\rho(t)\preceq \frac{t}{(\log t)^{1/p}}.$$
The following proposition answers negatively to this question. We
say that a function $f$ is sublinear if $f(t)/t\to 0$ when $t\to
\infty$.
\begin{PropIntro}\textnormal{(See Proposition~\ref{Cpbehavior})}
For any increasing sublinear function $h:\R_+\rightarrow \R_+$ and
every $1\leq p<\infty$, there exists a nondecreasing function $f$
satisfying $(C_p)$, a constant $c>0$ and a increasing sequence of
integers $(n_i)$ such that
$$ f(n_i)\geq ch(n_i), \quad \forall i\in \N.$$
\end{PropIntro}
In particular, it follows from Theorem~\ref{thm1} that the
compression $\rho$ of a uniform embedding of a 3-regular tree in a
Hilbert space does not satisfy any {\it a priori} majoration by
any sublinear function.

\subsection{Isoperimetry and compression}

To prove Theorem~\ref{thm1}, we observe a general relation between
the $L^p$-isoperimetry inside balls and the $L^p$-compression. Let
$G$ be a locally compact compactly generated group and consider
some compact symmetric generating subset $S$. For every $g\in G$,
write\footnote{We write $\tilde{\nabla}$ instead of $\nabla$
because this is not a``metric" gradient. The gradient associated
to the metric structure would be the right gradient: $|\nabla
\varphi|(g)=\sup_{s\in S}|\varphi(gs)-\varphi(g)|.$ This
distinction is only important when the group is non-unimodular.}
$$|\tilde{\nabla}\varphi|(g)=\sup_{s\in S}|\varphi(sg)-\varphi(g)|.$$
Let $2\leq p<\infty$ and let us call the $L^p$-isoperimetric
profile inside balls the nondecreasing function $J^b_{G,p}$
defined by
$$J^b_{G,p}(t)=\sup_{\varphi} \frac{\| \varphi\|_p}{\|\tilde{\nabla}\varphi\|_p},$$
where the supremum is taken over all measurable functions in
$L^p(G)$ with support in the ball $B(1,t)$. Note that the group
$G$ is amenable if and only if
$\lim_{t\to\infty}J^b_{G,p}(t)=\infty.$ Theorem~\ref{thm1} results
from the two following theorems.

\begin{ThmIntro}\label{isopthm}\textnormal{(see Theorem~\ref{Lthm})}
Let $G$ be a group of class $(\mathcal{L})$. Then
$J^b_{G,p}(t)\approx t$.
\end{ThmIntro}

\begin{ThmIntro}\label{CJpthm}\textnormal{(see Corollary~\ref{corCJp})}
Let $G$ be a locally compact compactly generated group and let $f$
be a nondecreasing function satisfying
\begin{equation}\tag{$CJ_p$}
\int_1^{\infty}\left(\frac
{f(t)}{J^b_{G,p}(t)}\right)^p\frac{dt}{t}<\infty
\end{equation}
for some $1< p<\infty$. Then there exists a $1$-cocycle $b\in
Z^1(G,\lambda_{G,p})$ whose compression $\rho$ satisfies
$\rho\succeq f$.
\end{ThmIntro}

Theorem~\ref{isopthm} may sound as a ``functional" property of
groups of class $(\mathcal{L})$. Nevertheless, our proof of this
result is based on a purely geometric construction. Namely, we
prove that these groups admit controlled F\o lner pairs (see
Definition \ref{Folnerpair}). In particular, when $p=1$ we obtain
the following corollary of Theorem~\ref{isopthm}, which has its
own interest.

\begin{ThmIntro}\textnormal{(See Remark
\ref{controlfoln} and  Theorem~\ref{Lthm})} Let $G$ be a group of
class $(\mathcal{L})$ and let $S$ be some compact generating
subset of $G$. Then $G$ admits a sequence of compact subsets
$(F_n)_{n\in \N}$ satisfying the two following conditions

\noindent(i) there is a constant $c>0$ such that
$$\mu(sF_n\vartriangle F_n)\leq c \mu(F_n)/n\quad \forall s\in S, \forall n\in \N;$$

\noindent(ii) for every $n\in \N$, $F_n$ is contained
\footnote{Actually, they also satisfy $S^{[cn]}\subset F_n$ for a
constant $c>0$.} in $S^n$.

\noindent In particular, $G$ admits a controlled F\o lner sequence
in the sense of \cite{coteva}.
\end{ThmIntro}
This theorem is a strengthening of the well-known construction by
Pittet \cite{Pittet}. It is stronger first because it does not
require the group to be unimodular, second because the control
(ii) of the diameter is really a new property that was not
satisfied in general by the sequences constructed in
\cite{Pittet}.

\subsection{Compression, subexponential growth, and random walks}
Let $\pi$ be a isometric $L^p$-representation of $G$. Denote by
$B_{\pi}(G)$ the supremum of all $\alpha$ such that there exists a
$1$-cocycle $b\in Z^1(G,\pi)$ whose compression $\rho$ satisfies
$\rho(t)\succeq t^{\alpha}$. Denote by $B_p(G)$ the supremum of
$B_{\pi}(G)$ over all isometric $L^p$-representations $\pi$. For
$p=2$, $B_2(G)=B(G)$ has been introduced in \cite{GK} where it was
called the equivariant Hilbert compression rate. On the other
hand, define
$$\alpha_{G,p}=\liminf_{t\to\infty}\frac{\log J^b_{G,p}(t)}{\log t}.$$
As a corollary of Theorem~\ref{thm1}, we have
\begin{CorIntro}\label{B(G)=1}
For every $1\leq p<\infty$, and every group $G$ of the class
$(\mathcal{L})$, we have $B_p(G)=1$.
\end{CorIntro}

The following result is a corollary of Theorem~\ref{CJpthm}.
\begin{CorIntro}\label{corintroCJp}\textnormal{(see Corollary~\ref{corCJp})}
Let $G$ be a locally compact compactly generated group. For every
$0< p<\infty$, we have
$$B_{\lambda_{G,p}}(G)\geq \alpha_{G,p}.$$
\end{CorIntro}
The interest of this corollary is illustrated by the two following
propositions. Recall the volume growth of $G$ is the $\approx$
equivalence class $V_G$ of the function $r\mapsto \mu(B(1,r))$.

\begin{PropIntro}\textnormal{(see Proposition~\ref{croissanceinterm})}
Assume that there exists $\beta< 1$ such that $V_G(r)\preceq
e^{r^{\beta}}$. Then $$\alpha_{G,p}\geq 1-\beta.$$
\end{PropIntro}

As an example we obtain that $B(G)\geq 0,19$ for the first
Grigorchuk's group (see \cite{Bart} for the best known upper bound
of the growth function of this group).

Let $G$ be a finitely generated group and let $\nu$ be a symmetric
finitely supported probability measure on $G$. Write
$\nu^{(n)}=\nu\ast\ldots\ast\nu$ ($n$ times). Recall that
$\nu^{(n)}(1)$ is the probability of return of the random walk
starting at $1$ whose probability transition is given by $\nu$.

\begin{PropIntro}\label{randomprop}\textnormal{(see Proposition~\ref{marcheprop})}
Assume that there exists $\gamma< 1$ such that
$\nu^{(n)}(1)\succeq e^{-n^{\gamma}}$. Then
$$\alpha_{G,2}\geq (1-\gamma)/2.$$
\end{PropIntro}

In \cite{PitSal}, it is proved that if $G$ is a finitely generated
extension
$$1\to K\to G\to N\to 1$$ where $K$ is abelian and $N$ is
abelian with $\Q$-rank $d$. Then
$$\limsup_{n}\log(-\log( \nu^{(n)}(1)))\leq 1-2/(d+2)$$ for any symmetric
finitely supported probability on $G$.

\begin{CorIntro}\label{metabcor}
Assume that $G$ is a finitely generated extension $1\to K\to G\to
N\to 1$ where $K$ is abelian and $N$ is abelian with $\Q$-rank
$d$. Then $$B(G)\geq 1/(d+2).$$ In particular, $B(G)>0$ for any
finitely generated metabelian group~$G$.
\end{CorIntro}

\subsection{The case of $\Z\wr\Z$}\label{ZwrZSection}

Combining the construction of Theorem~\ref{thm1} for $C_2\wr\Z$
with the cocycle induced by the morphism of $\Z^{(\Z)}\to
\ell^p(\Z)$, we obtain (see Proposition~\ref{ZwrZProp} for the
details).

\begin{ThmIntro}\label{thmZwrZ}
Fix some $1\leq p<\infty$. Let $G=\Z\wr \Z$ and let $f$ be an
increasing function $f:\R_+\to \R_+$ satisfying
\begin{equation}\tag{$C_p$}\int_1^{\infty}\left(\frac{f(t)}{t^{p/(2p-1)}}\right)^p\frac{dt}{t}<\infty.
\end{equation}
Then there exists a $1$-cocycle $b\in Z^1(G,\lambda_{G,p})$ whose
compression $\rho$ satisfies
$$\rho\succeq f.$$
In particular, $$B_p(\Z\wr\Z)\geq \frac{p}{2p-1}.$$
\end{ThmIntro}

In a previous version of this paper, we stated the lower bound
$B(\Z\wr\Z)\geq 2/3$, but the proof that we gave relied on a wrong
version of Proposition~\ref{randomprop} (we stated
$\alpha_{G,2}\geq 1-\gamma$, which is wrong as shown by a
counter-example in \cite{NP}). The mistake, together with a proof
of the full statement $B_p(\Z\wr\Z)\geq \frac{p}{2p-1}$ (see
\cite[Lemma~7.8]{NP}) was communicated to us by Naor and Peres.
The proof that we propose here is essentially the same as the one
of \cite{NP}. However, we would like to mention that we already
knew this proof\footnote{Actually, we knew it even before
proposing the wrong proof relying on Proposition~\ref{randomprop},
but we had chosen to include the latter here as it was a direct
(and nice) consequence of the methods of the present paper.} that
we dedicated to another paper with more general estimates on the
$L^p$-compression of wreath products. 

\subsection{Questions}

\begin{que}{\bf (Condition $(C_p)$ for nilpotent connected Lie groups.)}
Let $N$ be a simply connected non-abelian nilpotent Lie group and
let $\rho$ be the compression of a $1$-cocycle with values in some
$L^p$-space (resp. of a uniform embedding into some $L^p$-space)
for $2\leq p<\infty$. Does $\rho$ always satisfies Condition
$(C_p)$?
\end{que}
A positive answer would lead to the optimality of
Theorem~\ref{thm1}. On the contrary, one should wonder if it is
possible, for any increasing sublinear function $f$, to find a
$1$-cocycle (resp. a uniform embedding) in $L^p$ with compression
$\rho\succeq f$. This would also be optimal since we know
\cite{Pauls} that $N$ cannot quasi-isometrically embed into any
uniformly convex Banach space. Namely, the main theorem in
\cite{Pauls} states that such a group cannot quasi-isometrically
embed into any CAT(0)-space. So this only directly applies to
Hilbert spaces, but the key argument, consisting in a comparison
between the large scale behavior of geodesics (not exactly in the
original spaces but in tangent cones of ultra-products of them) is
still valid if the target space is a Banach space with unique
geodesics, a property satisfied by uniformly convex Banach spaces.

\begin{que}{\bf (Quasi-isometric embeddings into $L^1$-spaces.)}
Which connected Lie groups quasi-isometrically embed into some
$L^1$-space?
\end{que}
It is easy to quasi-isometrically embed a simplicial tree $T$ into
$\ell^1$ (see for instance \cite{GK}). In \cite{BuSc,BuSc'}, it is
proved that every semisimple Lie group of rank $1$
quasi-isometrically embeds into a finite product of simplicial
trees, hence into a $\ell^1$-space. The above question is of
particular interest for simply-connected non-abelian nilpotent Lie
groups since they do not quasi-isometrically embed into any finite
product of trees. Kleiner and Cheeger recently announced a proof
that the Heisenberg group cannot quasi-isometrically embed into
any $L^1$-space.

\begin{que}\label{Q4} If $G$ is an amenable group, is it true that
$$B_p(G)=\alpha_{G,p}?$$
\end{que}
We conjecture that this is true for $\Z\wr \Z$, i.e. that
$B(\Z\wr\Z)=2/3$. A first step to prove this is done by
Proposition~\ref{reguliere} which, applied to $G=\Z\wr\Z$ says
that
$$B(\Z\wr\Z)=B_{\lambda_{G,2}}(\Z\wr\Z).$$
As a variant of the above question, we may wonder if the weaker
equality $B_{\lambda_{G,p}}(G)=\alpha_{G,p}$ holds, in other words
if Corollary~\ref{corintroCJp} is optimal for all amenable groups.
Possible counterexamples would be wreath products of the form
$G=\Z\wr H$ where $H$ has non-linear growth (e.g. $H=\Z^2$).
\begin{que}
Does there exist an amenable group $G$ with $B(G)=0$?
\end{que}
A candidate would be the wreath product $\Z\wr (\Z\wr \Z)$ since
the probability of return of any non-degenerate random walk in
this group satisfies
$$\nu^{(n)}(1)\prec e^{-n^{\gamma}}$$
for every $\gamma<1$ \cite[Theorem~2]{Erschler'}. It is proved in
\cite{AGS} that $B(\Z\wr(\Z\wr\Z))\leq 1/2.$

\begin{que}\label{q6}
Let $G$ be a compactly generated locally compact group. If $G$
admits an isometric action on some $L^p$-space, $p\geq 2$, with
compression $\rho(t)\succ t^{1/p}$, does it imply that $G$ is
amenable?
\end{que}

Recall that this was proved in \cite{GK,coteva} for $p=2.$ The
generalization to every $p\geq 2$ would be of great interest. For
instance, this would prove the optimality of a recent result of Yu
\cite{Yu} saying that every finitely generated hyperbolic group
admits a proper isometric action on some $\ell^p$-space for large
$p$ enough, with\footnote{This is clear in the proof.} compression
$\rho(t)\approx t^{1/p}$.


\bigskip

\noindent \textbf{Acknowledgments.} First, I would like to thank
Assaf Naor and Yuval Peres for showing me their paper \cite{NP}
where they noticed a mistake in my earlier proof of
$B(\Z\wr\Z)\geq 2/3$, and proposed a correct proof of this result
(see Section~\ref{ZwrZSection}). Let us also mention that in their
paper, Naor and Peres answer to many of our questions, in
particular they answer negatively to Question~\ref{Q4}. On the
other hand, in \cite{ANP}, Austin, Naor and Peres prove our
conjecture that $B(\Z\wr\Z)=2/3$.

I am indebted to Yves de Cornulier for his critical reading of the
manuscrit and for numerous valuable discussions. I also thank
Pierre de la Harpe and Alain Valette for their useful remarks and
corrections. I am also grateful to Mark Sapir, Swiatoslaw Gal, and
Guillaume Aubrun for interesting discussions.

\section{Preliminaries}
\subsection{Compression}
Let us recall some definitions. Let $(X,d_X)$ and $(Y,d_Y)$ be
metric spaces. A map $F: X\to Y$ is called a uniform embedding of
$X$ into $Y$ if
$$d_X(x,y)\to \infty\quad \Longleftrightarrow \quad d_Y(F(x),F(y))\to
\infty.$$

Note that this property only concerns the large-scale geometry. A
metric space $(X,d)$ is called {\it quasi-geodesic} if there
exist
$\delta>0$ and $\gamma\geq 1$ such that for all $x,y\in X$, there
exists a chain $x=x_0,x_1,\ldots,x_n=y$ satisfying:
$$\sum_{k=1}^n d(x_{k-1},x_k)\leq \gamma d(x,y),$$
$$\forall k=1,\ldots,n,\quad d(x_{k-1},x_k)\leq \delta.$$
If $X$ is quasi-geodesic and if $F: X\to Y$ is a uniform
embedding, then it is easy to see that $F$ is large-scale
Lipschitz, i.e. there exists $C\geq 1$ such that
$$\forall x,y\in X,\quad d_Y(F(x),F(y))\leq Cd_X(x,y)+C.$$
Nevertheless, such a map is not necessarily large scale
bi-Lipschitz (in other words, quasi-isometric).
\begin{defn}\label{compressiondef}
We define the compression $\rho:\R_+\to [0,\infty]$ of a map
$F:X\to Y$ by
$$\forall t>0, \quad \rho(t)=\inf_{d_X(x,y)\geq
t}d_Y(F(x),F(y)).$$
\end{defn}
Clearly, if $F$ is large-scale Lipschitz, then $\rho(t)\preceq t$.

\subsection{Length functions on a group}
Now, let $G$ be a group. A length function on $G$ is a function
$L:G\to \R_+$ satisfying $L(1)=0$, $L(gh)\leq L(g)+L(h),$ and
$L(g)=L(g^{-1}).$ If $L$ is a length function, then
$d(g,h)=L(g^{-1}h)$ defines a left-invariant pseudo-metric on $G$.
Conversely, if $d$ is a left-invariant pseudo-metric on $G$, then
$L(g)=d(1,g)$ defines a length function on $G$.

Let $G$ be a locally compact compactly generated group and let $S$
be some compact symmetric generating subset of $G$. Equip $G$ with
a proper, quasi-geodesic length function by
$$|g|_S=\inf\{n\in \N: g\in S^n\}.$$
Denote $d_S$ the associated left-invariant distance. Note that any
proper, quasi-geodesic left-invariant metric is quasi-isometric to
$d_S$, and so belongs to the same ``asymptotic class".

\subsection{Affine isometric actions and first cohomology}

Let $G$ be a locally compact group, and $\pi$ a isometric
representation (always assumed continuous) on a Banach space
$E=E_{\pi}$. The space $Z^1(G,\pi)$  is defined as the set of
continuous functions $b:G\to E$ satisfying, for all $g,h$ in $G$,
the 1-cocycle condition $b(gh)=\pi(g)b(h)+b(g)$. Observe that,
given a continuous function $b:G\to\mathcal{H}$, the condition
$b\in Z^1(G,\pi)$ is equivalent to saying that $G$ acts by affine
isometries on $\mathcal{H}$ by $\alpha(g)v=\pi(g)v+b(g)$. The
space $Z^1(G,\pi)$ is endowed with the topology of uniform
convergence on compact subsets.

The subspace of coboundaries $B^1(G,\pi)$ is the subspace (not
necessarily closed) of $Z^1(G,\pi)$ consisting of functions of the
form $g\mapsto v-\pi(g)v$ for some $v\in E$. In terms of affine
actions, $B^1(G,\pi)$ is the subspace of affine actions fixing a
point.

The first cohomology space of $\pi$ is defined as the quotient
space $$H^1(G,\pi)=Z^1(G,\pi)/B^1(G,\pi).$$

Note that if $b\in Z^1(G,\pi)$, the map $(g,h)\to\| b(g)-b(h)\|$
defines a left-invariant pseudo-distance on $G$. Therefore the
compression of a $1$-cocycle $b: (G,d_S)\to E$ is simply given by
$$\rho(t)=\inf_{|g|_S\geq t}\| b(g)\|.$$
The compression of an affine isometric action is defined as the
compression of the corresponding $1$-cocycle.
\begin{rem}
When the space $E$ is a Hilbert space\footnote{The same proof
holds for uniformly convex Banach spaces.}, it is well-known
\cite[§4.a]{HV} that $b\in B^1(G,\pi)$ if and only if $b$ is
bounded on $G$.
\end{rem}

\section{The maximal $L^p$-compression functions $M\rho_{G,p}$ and $M\rho_{\lambda_{G,p}}$}\label{W}

\subsection{Definitions and general results}
Let $(G,d_S,\mu)$ be a locally compact compactly generated group,
generated by some compact symmetric subset $S$ and equipped with a
left Haar measure $\mu$. Denote by $Z^1(G,p)$ the collection of
all $1$-cocycles with values in any $L^p$-representation of $G$.
Denote by $\rho_b$ the compression function of a $1$-cocycle $b\in
Z^1(G,p)$.
\begin{defn}\label{maxcomp} We call maximal $L^p$-compression function of $G$
the nondecreasing function $M\rho_{G,p}$  defined by
$$M\rho_{G,p}(t)=\sup\left\{\rho_b(t):\;b\in Z^{1}(G,p),\; \sup_{s\in S}\|b(s)\|\leq 1\right\}.$$
We call maximal regular $L^p$-compression function of $G$ the
nondecreasing function $M\rho_{\lambda_{G,p}}$ defined by
$$M\rho_{\lambda_{G,p}}=\sup\left\{\rho_b(t):\;b\in Z^{1}(G,\lambda_{G,p}),\; \sup_{s\in S}\|b(s)\|\leq 1\right\}.$$
\end{defn}
Note that the asymptotic behaviors of both $M\rho_{G,p}$ and
$M\rho_{\lambda_{G,p}}$ do not depend on the choice of the compact
generating set $S$. Moreover, we have $$
M\rho_{\lambda_{G,p}}(t)\leq M\rho_{G,p}(t)\preceq t.$$

Let $\varphi$ be a measurable function on $G$ such that
$\varphi-\lambda(s)\varphi\in L^p(G)$ for every $s\in S$. For
every $t>0$, define
$$\Var_p(\varphi,t)=\inf_{|g|_S\geq t}\| \varphi-\lambda(g)\varphi\|_p$$
The function $\varphi$ and $p$ being fixed, the map $t\mapsto
\Var_p(\varphi,t)$ is nondecreasing.
\begin{prop}
We have
$$M\rho_{\lambda_{G,p}}(t)=\sup_{\|\tilde{\nabla}\varphi\|_p\leq 1}\Var_p(\varphi,t).$$
\end{prop}

\bpr We trivially have
$$M\rho_{\lambda_{G,p}}(t)\geq \sup_{\|\tilde{\nabla}\varphi\|_p\leq
1}\Var_p(\varphi,t).$$ Let $b$ be an element of
$Z^{1}(G,\lambda_{G,p})$. By a standard argument of
convolution\footnote{One can convolute $b(g)$, for every $g$, on
the right by a Dirac approximation.}, one can approximate $b$ by a
cocycle $b'$ such that $x\to b'(g)(x)$ is continuous for every $g$
in $G$. Hence, we can assume that $b(g)$ is continuous for every
$g$ in $G$. Now, setting $\varphi(g)=b(g)(g)$, we define a
measurable function satisfying
$$b(g)=\varphi-\lambda(g)\varphi.$$
So we have $$\rho(t)=\Var_p(\varphi,t)\leq
M\rho_{\lambda_{G,p}}(t)$$ where $\rho$ is the compression of $b$.
\epr

\

\begin{rem}
It is not difficult to prove that the asymptotic behavior of
$M\rho_{\lambda_{G,p}}$ is invariant under quasi-isometry between
finitely generated groups. 
\end{rem}

\begin{prop}
The group $G$ admits a proper\footnote{For $p=2$, this means that
$G$ is a-T-menable if and only if $M\rho_{G,2}$ goes to infinity.
It should be compared to the role played by the H-metric (see
\cite[§~2.6]{Corn} and §~7.4) for Property (T).} $1$-cocycle with
values in some $L^p$-representation if and only if
$M\rho_{G,p}(t)$ goes to infinity when $t\to\infty.$
\end{prop}

\bpr The ``only if" part is trivial. Assume that $M\rho_{G,p}(t)$
goes to infinity. Let $(t_k)$ be an increasing sequence growing
fast enough so that
$$\sum_{k\in \N}\frac{1}{t_k^p}<\infty.$$ For every $k\in \N$,
choose some $b_k\in Z^1(G,p)$ whose compression $\rho_k$ satisfies
$$\rho_k(t_k)\geq \frac{M\rho_{G,p}(t_k)}{2}$$ and such that
$$\sup_{s\in S}\|b_k(s)\|\leq 1.$$
Clearly, we can define a $1$-cocycle $b\in Z^1(G,p)$ by
$$b=\oplus_{k}^{\ell^p}\frac{1}{t_k}b_k.$$
That is, if for every $k$, $b_k$ takes values in the
representation $\pi_k$, then $b$ takes values in the direct sum
$\oplus_k^{\ell^p} \pi_k$. Now, observe that for $|g|\geq t_k$ and
$j\leq k$, we have $\|b_j(g)\|\geq 1/2$, so that $$\|b(g)\|^p\geq
k/2^p.$$ Thus the cocycle $b$ is proper.\epr

\

The following proposition, which is a quantitative version of the
previous one, plays a crucial role in the sequel.
\begin{prop}\label{Wp/compression}
Let $f:\R_+\to \R_+$ be a nondecreasing map satisfying
\begin{equation}\tag{$CM_p$}
\int_1^{\infty}\left(\frac{f(t)}{M\rho_{G,p}(t)}\right)^p\frac{dt}{t}<\infty,
\end{equation}
Then,

\noindent(1) there exists a $1$-cocycle $b\in Z^1(G,p)$ such that
$$\rho\succeq f;$$

\noindent(2) if one replace $M\rho_{G,p}$ by
$M\rho_{\lambda_{G,p}}$ in Condition~$(CM_p)$, then $b$ can be
chosen in  $Z^1(G,\lambda_{G,p})$.

\end{prop}
\noindent{\bf Proof of (1):} For every $k\in \N$, choose some
$b_k\in Z^1(G,p)$ (for (2), we take $b_k\in Z^1(G,\lambda_{G,p})$)
whose compression $\rho_k$ satisfies
$$\rho_k(2^{k+1})\geq \frac{M\rho_{G,p}(2^{k+1})}{2}$$ and such that
$$\sup_{s\in S}\|b_k(s)\|\leq 1.$$
Then define another sequence of cocycles $\tilde{b}_k\in Z^1(G,p)$
by
$$\tilde{b}_k=\frac{f(2^k)}{M\rho_{G,p}(2^{k+1})}b_k.$$
Since $M\rho_{G,p}$ and $f$ are nondecreasing, for any $2^k\leq
t\leq 2^{k+1}$, we have
$$\frac{f(2^k)}{M\rho_{G,p}(2^{k+1})}\leq \frac{f(t)}{M\rho_{G,p}(t)}.$$
Hence, for $s\in S$,
\begin{eqnarray*}
\sum_{k}\| \tilde{b}_k(s)\|_p^p &\leq & \sum_{k}
\left(\frac{f(2^k)}{M\rho_{G,p}(2^{k+1})}\right)^p\\
                                    & \leq &
2\int_{1}^{\infty}\left(\frac{f(t)}{M\rho_{G,p}(t)}\right)^p\frac{dt}{t}<\infty
\end{eqnarray*}
So we can define a $1$-cocycle on $b\in Z^1(G,p)$ by
\begin{equation}\label{cocycle}
b=\oplus_{k}\tilde{b}_k.
\end{equation}

On the other hand, if $|g|_S\geq 2^{k+1}$, then
\begin{eqnarray*}
\|
b(g)\|_p & \geq & \| \tilde{b}_k(g)\|_p\\
                & \geq &
\frac{f(2^k)}{M\rho_{\lambda_{G,p}}(2^{k+1})}\rho_k(2^{k+1})\\
                & \geq & f(2^k)
\end{eqnarray*}
So if $\rho$ is the compression of the $1$-cocycle $b$, we have
$\rho\succeq f$.

\

\noindent{\bf Proof of (2):} We keep the previous notation. Assume
that $f$ satisfies
$$\int_1^{\infty}\left(\frac{f(t)}{M\rho_{\lambda_{G,p}}(t)}\right)^p\frac{dt}{t}<\infty.$$
The cocycle $b$ provided by the proof of (1) has the expected
compression but it takes values in an infinite direct sum of
regular representation $\lambda_{G,p}$. Now, we would like to
replace the direct sum $b=\oplus_{k}b_k$ by a mere sum, in order
to obtain a cocycle in $Z^1(G,\lambda_{G,p})$. Since $G$ is not
assumed unimodular, the measure $\mu$ is not necessarily
right-invariant. However, one can define a isometric
representation $r_{G,p}$ on $L^p(G)$, called the right regular
representation by
$$r_{G,p}(g)\varphi=\Delta(g)^{-1}\varphi(\cdot g)\quad \forall \varphi\in
L^p(G),$$where $\Delta$ is the modular function of $G$. We will
use the following well-known property of the representation
$r_{G,p}$, for $p>1$. To simplify, let us write $r(g)$ instead of
$r_{G,p}(g).$ For every $(\varphi,\psi)\in L^p(G)\times L^p(G)$,
we have
\begin{equation}\label{C0}
\lim_{|g|\to\infty}\|
r(g)\varphi+\psi\|_p^p=\|\varphi\|_p^p+\|\psi\|_p^p.
\end{equation}
Moreover, this limit is uniform on compact subsets of
$(L^p(G))^2$. As $r_{G,p}$ and $\lambda_{G,p}$ commute, $r_{G,p}$
acts by isometries on $Z^1(G,\lambda_{G,p})$.
\begin{lem}\label{lem}
There exists a sequence $(g_k)$ of elements of $G$ such that
$b'=\sum r(g_k)b_k$ defines a cocycle in $Z^1(G,\lambda_{G,p})$
and such that
\begin{equation}\label{eq1}
\left| \|b'(g)\|_p^p-\left\|
\sum_{j=0}^{k-1}r(g_j)b_j(g)\right\|_p^p -\sum_{j\geq k} \|
b_j(g)\|_p^p\right|\leq 1
\end{equation}
for any $k$ large enough and every $g\in B(1,2^{k+2})$.
\end{lem}
\noindent{\bf Proof of Lemma \ref{lem}.} By an immediate
induction, using (\ref{C0}), we construct a sequence $(g_k)\in
G^{\N}$ satisfying, for every $K\geq 0,s\in S,$
$$\| \sum_{k=0}^K r(g_k)b_k(s)\|_p^p\leq
\sum_{k=0}^K\| b_k(s)\|_p^p+\sum_{k=0}^K 2^{-k-1}\leq 1,$$ which
implies that $b'$ is a well-defined $1$-cocycle in
$Z^1(G,\lambda_{G,p})$. Similarly, one can choose $(g_k)$
satisfying the additional property that, for every $k\in \N,$
$|g|\leq 2^{k+2}$,
$$\left|\|\sum_{j=0}^k
r(g_j)b_j(g)\|_p^p-\|\sum_{j=0}^{k-1} r(g_j) b_j(g)\|_p^p-\|
b_k(g)\|_p^p\right|\leq 2^{-k-1}.$$Fixing $k\in \N$, an immediate
induction over $K$ shows that for every $|g|\leq 2^{k+2}$ and
every $K\geq k$,
$$\left|\|\sum_{j=0}^K r(g_j)b_j(g)\|_p^p-\|
\sum_{j=0}^{k-1}r(g_j)b_j(g)\|_p^p -\sum_{j=k}^K \|
b_j(g)\|_p^p\right|\leq \sum_{j=k}^K 2^{-j-1}.$$ This proves
(\ref{eq1}).\epr

\

By the lemma, for $|g|\leq 2^{k+2}$,
\begin{eqnarray*}
\| b'(g)\|_p^p & \geq & \| b_k(g) \|_p^p-1.
\end{eqnarray*}
Then, for $2^{k+1}\leq |g|\leq 2^{k+2}$, we have
\begin{eqnarray*}
\| b'(g)\|_p^p & \geq & f(2^k)-1
\end{eqnarray*}
Therefore, the compression $\rho'$ of $b'$ satisfies
$$\rho'\succeq f$$
and we are done.\epr

\

We have the following immediate consequence.
\begin{cor}
For every $1\leq p<\infty,$
$$B(G,p)= \liminf_{t\to \infty}\frac{\log M\rho_{G,p}(t)}{\log
t}.$$
\end{cor}

\begin{ex}
Let $F_r$ be the free group of rank $r\geq 2$ and let $A(F_r)$ be
the set of edges of the Cayley graph of $F_r$ associated to the
standard set of generators. The standard isometric affine action
of $F_r$ on $\ell^p(A(F_r))$, whose linear part is isomorphic to a
direct sum $\lambda_{G,p}\oplus_{\ell^p}\ldots\oplus_{\ell^p}
\lambda_{G,p}$ of $r$ copies of $\lambda_{G,p}$ has compression
$\approx t$. This shows that $M\rho_{\lambda_{F_r,p}}(t)\succeq
t^{1/p}$.
\end{ex}

\subsection{Reduction to the regular representation for $p=2$}
In the Hilbert case, we prove that if a group admits a $1$-cocycle
with large enough compression, then
$M\rho_{G,2}=M\rho_{\lambda_{G,2}}$. This result is mainly
motivated by Question~\ref{Q4} since it implies that
$$B(\Z\wr\Z)=B_{\lambda_{G,2}}(\Z\wr\Z).$$

\begin{prop}\label{reguliere}
Let $\pi$ be a unitary representation of the group $G$ on a
Hilbert space $\HH$ and let $b\in Z^1(G,\pi)$ be a cocycle whose
compression $\rho$ satisfies
$$\rho(t)\succ t^{1/2}.$$
Then\footnote{Note that the hypotheses of the proposition also
imply that $G$ is amenable \cite[Theorem~4.1]{coteva,GK}.},
$$\rho\preceq M\rho_{\lambda_{G,2}}.$$ In particular,
$$M\rho_{2}=M\rho_{\lambda_{G,2}}.$$
\end{prop}
combining with Proposition~\ref{Wp/compression}, we obtain
\begin{cor}
With the same hypotheses, we have
$$B(G)=B(G,\lambda_{G,2})=\liminf_{t\to \infty}\frac{\log
M\rho_{\lambda_{G,2}}(t)}{\log t}.$$
\end{cor}

\noindent{\bf Proof of Proposition~\ref{reguliere}.} For every
$t>0$, define
$$\varphi_t(g)=e^{-\| b(g)\|^{2}/t^2}.$$
By Schoenberg's Theorem \cite[Appendix C]{BHV}, $\varphi_t$ is
positive definite. It is easy to prove that $\varphi_t$ is
square-summable (see \cite[Theorem~4.1]{coteva}). By
\cite[Théorème 13.8.6]{Dix}, it follows that there exists a
positive definite, square-summable function $\psi_t$ on $G$ such
that $\varphi_t=\psi_t\ast\psi_t$, where $\ast$ denotes the
convolution product. In other words, $\varphi_t=\langle \lambda(g)
\psi_t,\psi_t\rangle.$ In particular,
$$\varphi_t(1)=1=\|\psi_t\|_2^2$$
and for every $s\in S$,
\begin{eqnarray*}
\|\psi_t-\lambda(s)\psi_t\|_2^2 & = & 2(\|\psi_t\|_2^2-\langle
\lambda(s)
\psi_t,\psi_t\rangle)\\
 & = & 2(1-\varphi_t(s))\\
 & = & 2(1-e^{-\| b(s)\|^2/t^2})\\
 & \preceq 1/t^2
\end{eqnarray*}
On the other hand, for $g$ such that $\rho(|g|_S)\geq t$, we have
\begin{eqnarray*}
\|\psi_t-\lambda(g)\psi_t\|_2^2
 & = & 2(1-e^{-\| b(g)\|^2/t^2})\\
 & \geq & 2(1-e^{-\rho(|g|_S)^2/t^2})\\
 & \geq & 2(1-1/e)
\end{eqnarray*}
So, we have
$$\frac{\|\psi_t-\lambda(g)\psi_t\|_2}{\| \tilde{\nabla}\psi_t\|_2}\geq ct$$
where $c$ is a constant. In other words,
$$\Var_{2}(\psi_t,\rho^{-1}(t))\geq ct.$$
It follows from the definitions that $M\rho_{\lambda_{G,2}}\succeq
\rho.$\epr

\section{$L^p$-isoperimetry inside balls}

\subsection{Comparing $J^b_{G,p}$ and $M\rho_{\lambda_{G,p}}$}

Let $G$ be a locally compact compactly generated group and let $S$
be a compact symmetric generating subset of $G$. Let $A$ be a
subset of the group $G$. One defines the $L^p$-isoperimetric
profile inside $A$ by
$$J_{p}(A)=\sup_{\varphi} \frac{\| \varphi\|_p}{\|\tilde{\nabla}\varphi\|_p}$$
where the supremum is taken over nonzero functions in $L^p(G)$
with support included in $A$.
\begin{defn}
The $L^p$-isoperimetric profile inside balls is the nondecreasing
function $J^b_{G,p}$ defined by
$$J^b_{G,p}(t)=J_p(B(1,t)).$$
\end{defn}
\begin{rem}
The usual $L^p$-isoperimetric profile of $G$ (see for example
\cite{Coulhon}) is defined by
$$j_{G,p}(t)=\sup_{\mu(A)=t}J_p(A).$$
Note that our notion of isoperimetric profile depends on the
diameter of the subsets instead of their measure.
\end{rem}
\begin{rem}
The asymptotic behavior of $J^b_{p,G}$ is invariant under
quasi-isometry between compactly generated groups \cite{tess"}. In
particular, it is also invariant under passing to a cocompact
lattice \cite{CS}.
\end{rem}
\begin{rem}
Using basic $L^p$-calculus, one can easily prove \cite{Coulhon}
that if $p\leq q$, then
$$(J^b_{G,p})^{p/q}\preceq J^b_{G,q}\preceq J^b_{G,p}.$$
\end{rem}

Now let us compare $J^b_{p,G}$ and $M\rho_{\lambda_{G,p}}$
introduced in §\;\ref{W}.
\begin{prop}\label{propWJ}
For every $2\leq p<\infty$, we have $$M\rho_{\lambda_{G,p}}\succeq
J^b_{G,p}.$$
\end{prop}
\bpr Fix some $t>0$ and choose some $\varphi\in L^p(X)$ whose
support lies in $B(1,t)$ such that
$$\frac{\|
\varphi\|_p}{\|\tilde{\nabla}\varphi\|_p}\geq J^b_{G,p}(t)/2.$$
Take $g\in G$ satisfying $|g|_S\geq 3t$. Note that $B(1,t)\cap
\lambda(g)B(1,t)=\emptyset.$ So $\varphi$ and $\lambda(g)\varphi$
have disjoint supports. In particular,
$$\|\varphi-\lambda(g)\varphi\|_p\geq \|\varphi\|_p$$
and
$$\|\tilde{\nabla}(\varphi-\lambda(g)\varphi)\|_p=2^{1/p}\|\tilde{\nabla}\varphi\|_p$$
This clearly implies the proposition.\epr

\

Combining with Proposition~\ref{Wp/compression}, we obtain
\begin{cor}\label{corCJp}
Let $f:\R_+\to \R_+$ a nondecreasing map be satisfying
\begin{equation}\tag{$CJ_p$}
\int_1^{\infty}\left(\frac{f(t)}{J^b_{G,p}(t)}\right)^p\frac{dt}{t}<\infty
\end{equation}
for some $1\leq p<\infty.$ Then there exists a $1$-cocycle $b$ in
$Z^1(G,\lambda_{G,p})$ such that
$$\rho\succeq f.$$
\end{cor}

\begin{que}\label{questionWJ}
For which groups $G$ do we have $M\rho_{\lambda_{G,p}}\approx
J^b_{G,p}$?
\end{que}
We show that the question has positive answer for groups of class
$(\mathcal{L})$. On the contrary, note that the group $G$ is
nonamenable if and only if $J^b_{G,p}$ is bounded. But we have
seen in the previous section that for a free group of rank $\geq
2$, $M\rho_{\lambda_{G,p}}(t)\succeq t^{1/p}.$ More generally, the
answer to Question~\ref{questionWJ} is no for every nonamenable
group admitting a proper $1$-cocycle with values in the regular
representation. This question is therefore only interesting for
amenable groups.


\subsection{Sequences of controlled F\o lner pairs}

In this section, we give a method, adapted\footnote{In
\cite{CoulGri}, the authors are interested in estimating the
$L^2$-isoperimetric profile of a group.} from \cite{CoulGri} to
estimate $J^b_p$.

\begin{defn}\label{Folnerpair}
Let $G$ be a compactly generated, locally compact group equipped
with a left invariant Haar measure $\mu$. Let $\alpha=(\alpha_n)$
be a nondecreasing sequence of integers. A sequence of
$\alpha$-controlled F\o lner pairs of $G$ is a family
$(H_{n},H'_n)$ where $H_n$ and $H_n'$ are nonempty compact subsets
of $G$ satisfying for some constant $C>0$ the following
conditions:

\noindent(1) $S^{\alpha_n}H_{n}\subset H_{n}'$

\noindent(2) $\mu(H'_{n})\leq C\mu(H_{n})$;

\noindent(3) $H'_n\in B(1,Cn)$

If $\alpha_n\approx n$, we call $(H_n,H'_n)$ a controlled sequence
of F\o lner pairs.
\end{defn}

\begin{prop}\label{control/profile}
Assume that $G$ admits a sequence of $\alpha$-controlled F\o lner
pairs. Then $$J^b_{G,p}\succeq \alpha.$$
\end{prop}
\bpr For every $n\in \N$, consider the function $\varphi_n:G\to
\R_+$ defined by
$$\varphi_n(g)=\min\{k\in \N: g\in S^k(H'_n)^c\}$$
where $A^c=G\smallsetminus A.$ Clearly, $\varphi_n$ is supported
in $H'_n$. It is easy to check that
$$\| \tilde{\nabla}\varphi_n\|_p\leq (\mu(H'_n))^{1/p}$$
and that
$$\| \varphi_n\|_p\geq
\alpha_n(\mu(H_n))^{1/p}.$$ Hence by (2),
$$\| \varphi_n\|_p\geq C^{-1/p}\alpha_n\| \tilde{\nabla}\varphi_n\|_p,$$
so we are done.\epr

\begin{rem}\label{controlfoln}
Note that if $H$ and $H'$ are subsets of $G$ such that
$S^kH\subset H'$ and $\mu(H')\leq C\mu(H)$, then there exists by
pigeonhole principle an integer $0\leq j\leq k-1$ such that
$$\mu(\partial S^jH)= \mu(S^{j+1}H\smallsetminus S^jH)\leq
\frac{C}{k}\mu(S^jH).$$ So in particular if $(H_n,H'_n)$ is a
$\alpha$-controlled sequence of F\o lner pairs, then there exists
a F\o lner sequence $(K_n)$ such that $H_n\subset K_n\subset H'_n$
and $$\frac{\mu(\partial K_n)}{\mu(K_n)}\leq C/\alpha_n.$$
Moreover, if $\alpha_n\approx n$, then one obtains a controlled
F\o lner sequence in the sense of \cite[Definition~4.8]{coteva}.
\end{rem}

\section{Isoperimetry in balls for groups of class
$(\mathcal{L})$}\label{goodfol}

The purpose of this section is to prove the following theorem.

\begin{thm}\label{Lthm}
Let $G$ be a group belonging to the class $(\mathcal{L})$. Then,
$G$ admits a controlled sequence of F\o lner pairs. In particular,
$J^b_{G,p}(t)\approx t$.
\end{thm}

Note that Theorem~\ref{thm1} follows from Theorem~\ref{Lthm} and
Corollary~\ref{corCJp}.

\subsection{Wreath products $F\wr \Z$ }\label{Folwr}

Let $F$ be a finite group. Consider the wreath product $G=F\wr
\Z=\Z\ltimes F^{(\Z)}$, the group law being defined as
$(n,f)(m,g)=(n+m,\tau_mf+g)$ where $\tau_mf(x)=f(m+x)$. As a set,
$G$ is a Cartesian product $\Z\times U$ where $U$ is the direct
sum $F^{(\Z)}=\bigoplus_{n\in \Z}F_n$ of copies $F_n$ of $F$. The
set $S=S_F\cup S_{\Z}F_n$, where $S_F=F_{0}$ and
$S_{\Z}=\{-1,0,1\}$ is clearly a symmetric generating set for $G$.
Define
$$H_{n}=I_{n}\times U_{n}$$ and
$$H'_n=I_{2n}\times U_{n}$$
where $U_{n}=F^{[-2n,2n]}$ and $I_n=[-n,n].$

\noindent Let us prove that $(H_{n},H'_n)_n$ is a sequence of
controlled F\o lner pairs. We therefore have to show that

\noindent(1) $S^nH_{n}\subset H_{n}'$

\noindent(2) $|H'_{n}|\leq 2|H_{n}|$;

\noindent(3) there exists $C>0$ such that $H'_n\subset B(1,Cn)$

Property (2) is trivial. To prove (1) and (3), recall that the
length of an element of $g=(k,u)$ of $G$ equals
$L(\gamma)+\sum_{h\in \Z}|u(h)|_F$ where $L(\gamma)$ is the length
of a shortest path $\gamma$ from $0$ to $k$ in $\Z$ passing
through every element of the support of $u$ (see \cite[Theorem
1.2]{Par}). In particular,
$$|(u,k)|_S\leq 2L(\gamma)$$
Thus, if $g\in H_{n}$, then $L(\gamma)\leq 30n$. So (3) follows.
On the other hand, if $g=(k,u)\in S^n$, then
$$|k|_{\Z}\leq L(\gamma)\leq n$$
and
$$\Supp(u)\subset I_n.$$
So $H_ng\subset H_n'$.\epr

\begin{rem}
Note that the proof still works replacing $\Z$ by any group with
linear growth. On the other hand, replacing it by a group of
polynomial growth of degree $d$ yields a sequence of
$n^{1/d}$-controlled F\o lner pairs. For instance, as a corollary,
we obtain that $B(F\wr\Z^d)\geq 1/d$.
\end{rem}

\subsection{Semidirect products
$\left(\R\oplus\bigoplus_{p\in
P}\Q_{p}\right)\rtimes_{\frac{m}{n}}\Z$.}\label{Folsemi}

Note that discrete groups of type (2) of the class $(\mathcal{L})$
are cocompact lattices of a group of the form
$$G=\Z\ltimes_{\frac{m}{n}}\left(\R\oplus\bigoplus_{p\in
P}\Q_{p}\right)$$ with $m,n$ coprime integers and $P$ a finite set
of primes (possibly infinite) dividing $mn$. To simplify notation,
we will only consider the case when $P=\{p\}$ is reduced to one
single prime, the generalization presenting no difficulty. The
case where $p=\infty$ will result from the case of connected Lie
groups (see next section) since $\Z\ltimes_{\frac{m}{n}}\R $
embeds as a closed cocompact subgroup of the group of positive
affine transformations $\R\ltimes\R$.

So consider the group $G=\Z\ltimes_{1/p}\Q_{p}$. Define a compact
symmetric generating set by $S=S_{\Q_p}\cup S_{\Z}$ where
$S_{\Q_p}=\Z_p$ and $S_{Z}=\{-1,0,1\}$. Define $(H_k,H'_k)$ by
$$H_{k}=I_{k}\times p^{-2k}\Z_p$$
and
$$H'_{k}=I_{2k}\times p^{-2k}\Z_p,$$
where $I_k=[-k, k]$. Using the same kind of arguments as
previously for $F\wr \Z$, one can prove easily that $(H_k,H'_k)$
is a controlled sequence of F\o lner pairs. \epr

\subsection{Amenable connected Lie groups}

Let $G$ be a solvable simply connected Lie group. Let $S$ be a
compact symmetric generating subset. In \cite{Osin}, it is proved
that $G$ admits a maximal normal connected subgroup such that the
quotient of $G$ by this subgroup has polynomial growth. This
subgroup is called the exponential radical and is denoted
$\Exp(G)$. We have $\Exp(G)\subset N,$ where $N$ is the maximal
nilpotent normal subgroup of $G$. Let $T$ be a compact symmetric
generating subset of $\Exp(G).$ An element $g\in G$ is called
strictly exponentially distorted if the $S$-length of $g^n$ grows
as $\log |n|$. The subset of strictly exponentially distorted
elements of $G$ coincides with $\Exp(G)$. That is,
$$\Exp(G)=\{g\in G:\; |g^n|_S\approx \log n\}.$$
Moreover, $\Exp(G)$ is strictly exponentially distorted in $G$ in
the sense that there exists $\beta\geq 1$ such that for every
$h\in \Exp(G),$
\begin{equation}\label{expdist}
\beta^{-1}\log(|h|_{T}+1)-\beta\leq |h|_S \leq
\beta\log(|h|_{T}+1)+\beta
\end{equation}
where $T$ is a compact symmetric generating subset of $\Exp(G).$

We will need the following two lemmas.

\begin{lem}\label{sectionlem}
Let $G$ be a locally compact group. Let $H$ be a closed normal
subgroup. Let $\lambda$ and $\nu$ be respectively left Haar
measures of $H$ and $G/H$. Let $i$ be a measurable left-section of
the projection $\pi:~G\to G/H$, i.e. $G=\sqcup_{x\in G/H}i(x)H$.
Identify $G$ with the cartesian product $G/H\times H$ via the map
$(x,h)\mapsto i(x)h.$ Then the product measure $\nu\otimes\lambda$
is a left Haar measure on $G$.
\end{lem}
\bpr We have to prove that $\nu\otimes\lambda$ is left-invariant
on $G$. Fix $g$ in $G$. Define a measurable map $\sigma_g$ from
$G/H$ to $H$ by
$$\sigma_g(x)=(i(\pi(g)x)^{-1}gi(x).$$
In other words, $\sigma_g(x)$ is the unique element of $H$ such
that
$$gi(x)=i(\pi(g)x)\sigma_g(x).$$
Let $\varphi~:G\to \R$ be a continuous, compactly supported
function. We have
$$\int_{G/H\times H}\varphi[gi(x)h]d\nu(x)d\lambda(h) =
\int_{G/H\times
H}\varphi[i(\pi(g)x)\sigma_g(x)h]d\nu(x)d\lambda(h)$$ As $\nu$ and
$\lambda$ are respectively left Haar measures on $G/H$ and $H$,
the Jacobian of the transformation
$(x,h)\mapsto(\pi(g)x,\sigma_g(x)h)$ is equal to $1$. Hence,
\begin{eqnarray*}
\int_{G/H\times
H}\varphi[i(\pi(g)x)\sigma_g(x)h]d\nu(x)d\lambda(h) &=&
\int_{G/H\times H}\varphi[i(x)h]d\nu(x)d\lambda(h).
\end{eqnarray*}
Thus $\nu\otimes\lambda$ is left-invariant.\epr

\begin{lem}\label{sectionlem'}
Let $G$ be a connected Lie group and $H$ be a normal subgroup.
Consider the projection $\pi:G\to G/H$. There exists a compact
generating set $S$ of $G$ and a $\sigma$-compact cross-section
$\sigma$ of $G/H$ inside $G$ such that $\sigma(\pi(S)^n)\subset
S^{n+1}$.
\end{lem}
\bpr Since $\pi$ is a submersion, there exists a compact
neighborhood $S$ of $1$ in $G$ such that $\pi(S)$ admits a
continuous cross-section $\sigma_1$ in $S$. Now, let $X$ be a
minimal (discrete) subset of $G/H$ satisfying $G/H=\cup_{x\in X}
x\pi(S)$. Since this covering is locally finite and $\pi(S)$ is
compact, one can construct by induction a partition $(A_x)_{x\in
X}$ of $G/H$ such that every $A_x$ is a constructible, and
therefore $\sigma$-compact subset of $x\pi(S)$. Let $\sigma_2:
X\to G$ be a cross-section of $X$ satisfying $\sigma_2(X\cap
\pi(S)^n)\subset S^n$. Now, for every $z\in A_x$, define
$$\sigma(z)=\sigma_2(x)\sigma_1(x^{-1}z).$$
Clearly, $\sigma$ satisfies to the hypotheses of the lemma.\epr

\

Equip the group $P=G/\Exp(G)$ with a Haar measure $\nu$ and with
the symmetric generating subset $\pi(S)$, where $\pi$ is the
projection on $P$. Assume that $S$ satisfies to the hypotheses of
Lemma~\ref{sectionlem'} and let $\sigma$ be a $\sigma$-compact
cross-section of $P$ inside $G$ such that
$\sigma(\pi(S)^{n})\subset S^{n+1}$. For every $n\in \N$, write
$F_n=\sigma(\pi(S)^{n})$. Let $\alpha$ be some large enough
positive number that we will determine later. Denote by $\lfloor
x\rfloor$ the integer part of a real number $x$. Define, for every
$n\in \N$,
$$H_n=S^nT^{\lfloor\exp(\alpha n)\rfloor}$$
and
$$H'_n=S^{2n}T^{\lfloor\exp(\alpha n)\rfloor}.$$
Note that $H'_n=S^nH_n$. On the other hand, since $\Exp(G)$ is
strictly exponentially distorted, there exists $a\geq 1$ only
depending on $\alpha$ and $\beta$ such that, for every $n\in \N$,
$$S^nT^{\lfloor\exp(\alpha n)\rfloor}\subset S^{an}.$$
Hence, to prove that $(H_n,H'_n)$ is a sequence of controlled F\o
lner pairs, it suffices to show that $\mu(H_n')\leq C\mu(H_n)$.
Consider another sequence $(A_n,A'_n)$ defined by, for every $n\in
\N^*$,
$$A_n=F_{n-1}T^{\lfloor\exp(\alpha n)\rfloor}$$
and
$$A'_n=F_{2n}T^{2\lfloor\exp(\alpha n)\rfloor}.$$
As $F_n$ is $\sigma$-compact, $A_n$ and $A'_n$ are measurable. To
compute the measures of $A_n$ and $A'_n$, we choose a
normalization of the Haar measure $\lambda$ on $\Exp(G)$ such that
the measure $\mu$ disintegrates over $\lambda$ and the pull-back
measure of $\nu$ on $\sigma(P)$ as in Lemma~\ref{sectionlem}. We
therefore obtain
$$\mu(A_n)=\nu(\pi(S)^{n-1})\lambda(T^{\lfloor\exp(\alpha n)\rfloor})$$
and
$$\mu(A'_n)=\nu(\pi(S)^{2n})\lambda(T^{2\lfloor\exp(\alpha n)\rfloor}).$$
Since $P$ and $\Exp(G)$ have both polynomial growth, there is a
constant $C$ such that, for every $n\in \N^*$,
$$\mu(A'_n)\leq C\mu(A_n).$$ So now, it suffices to prove that
$$A_n\subset H_n\subset H'_n\subset A'_n,$$
where the only nontrivial inclusion is $H_n'\subset A'_n$. Let
$g\in S^{2n}$; let $f\in F_{2n}$ be such that $\pi(g)=\pi(f).$
Since $F_{2n}\subset S^{2n+2}\subset S^{3n},$ $$gf^{-1}\in
S^{6n}\cap \Exp(G).$$ On the other hand, by (\ref{expdist}),
$$S^{6n}\cap \Exp(G)\subset T^{2\lfloor\exp(6\beta n)\rfloor}.$$
Therefore, for every $n\in \N^*$,
$$H'_n\subset F_{2n}T^{2\lfloor\exp(6\beta n)\rfloor}T^{\lfloor\exp(\alpha n)\rfloor}=F_{2n}T^{2\lfloor\exp(6\beta n)\rfloor+\lfloor\exp(\alpha n)\rfloor}.$$
Hence, choosing $\alpha\geq 6\beta+\log 2$, we have
$$H'_n\subset F_{2n}T^{2\lfloor\exp(\alpha n)\rfloor}= A'_n,$$
and we are done.\epr

\section{On embedding of finite trees into uniformly convex Banach spaces}

\begin{defn}
A Banach space $X$ is called $q$-uniformly convex ($q>0$) if there
is a constant $a>0$ such that for any two points $x,y$ in the unit
sphere satisfying $\|x-y\|\geq \eps$, we have
$$\left\|\frac{x+y}{2}\right\|\leq 1-a\eps^q.$$
\end{defn}
Note that by a theorem of Pisier \cite{Pisier}, every uniformly
convex Banach space is isomorphic to some $q$-uniformly convex
Banach space.

In this section, we prove that the compression of a Lipschitz
embedding of a finite binary rooted tree into a $q$-uniformly
convex space $X$ always satisfies condition $(C_q).$
Theorem~\ref{finitetreethm} follows from the fact that a
$L^p$-space is $\max\{p,2\}$-uniformly convex.

\begin{thm}\label{embedtreethm}
Let $T_J$ be the binary rooted tree of depth $J$ and let $1<
q<\infty$. Let $F$ be a $1$-Lipschitz map from $T_J$ to some
$q$-uniformly convex Banach space $X$ and let $\rho$ be the
compression of $F$. Then there exists $C=C(q)<\infty$ such that
\begin{equation}\label{thmineq}
\int_{1}^{2J}\left(\frac{\rho(t)}{t}\right)^q\frac{dt}{t}\leq C.
\end{equation}
\end{thm}

\begin{cor}\label{Bourgainthm}
Let $F$ be any uniform embedding of the 3-regular tree $T$ into
some $q$-uniformly convex Banach space. Then the compression
$\rho$ of $F$ satisfies Condition $(C_q)$.\epr
\end{cor}

As a corollary, we also reobtain the theorem of Bourgain.
\begin{cor} \cite[Theorem 1]{Bourgain}
With the notation of Theorem~\ref{embedtreethm}, there exists at
least two vertices $x$ and $y$ in $T_J$ such that
$$\frac{\| F(x)-F(y)\|}{d(x,y)} \leq \left(\frac{C}{\log J}\right)^{1/q}.$$
\end{cor}
\bpr For every $1\leq t\leq 2J$, there exist $z,z'\in T_J$,
$d(z,z')\geq t$ such that:
$$\frac{\rho(t)}{t}=\frac{\| F(z)-F(z')\|}{t}\geq \frac{\|
F(z)-F(z')\|}{d(z,z')}.$$ Therefore
$$\min_{z\neq z'\in T_J}\frac{\|
F(z)-F(z')\|}{d(z,z')}\leq \min_{1\leq u\leq
2J}\frac{\rho(u)}{u}.$$ But by (\ref{thmineq})
$$\left(\min_{1\leq u\leq
2J}\frac{\rho(u)}{u}\right)^q\int_{1}^{2J}\frac{1}{t}dt\leq
\int_{1}^{2J}\left(\frac{\rho(t)}{t}\right)^q\frac{dt}{t}\leq C.$$
We then have
$$\min_{z\neq z'\in T_J}\frac{\|
F(z)-F(z')\|}{d(z,z')}\leq \left(\frac{C}{\log
J}\right)^{1/q}.\blacksquare
$$

\noindent{\bf Proof of Theorem~\ref{embedtreethm}.}  Since the
proof follows closely the proof of \cite[Theorem 1]{Bourgain}, we
keep the same notation to allow the reader to compare them. For
$j=1,2\ldots,$ denote $\Omega_j=\{-1,1\}^j$ and
$T_j=\bigcup_{j'\leq j}\Omega_{j'}$. Thus $T_j$ is the finite tree
with depth $j$. Denote $d$ the tree-distance on $T_j$.

\begin{lem}\cite[Proposition 2.4]{Pisier}
There exists $C=C(q)<\infty$ such that if $(\xi_{s})_{s\in \N}$ is
an $X$-valued martingale on some probability space $\Omega$, then
\begin{equation}\label{ineqlem1}
\sum_{s}\| \xi_{s+1}-\xi_{s}\|_q^q\leq C\sup_{s}\| \xi_{s}\|_q^q
\end{equation}
where $\|\;\|_q$ stands for the norm in $L^q_{X}(\Omega).$
\end{lem}\label{lem1}

Lemma~\ref{lem1} is used to prove
\begin{lem}
If $x_1,\ldots,x_J$, with $J=2^r$, is a finite system of vectors
in $X$, then
\begin{equation}\label{ineqlem2}
\sum_{s=1}^r2^{-qs}\min_{2^s\leq j\leq J-2^s} \|
2x_{j}-x_{j-2^s}-x_{j+2^s}\|^q\leq C\sup_{1\leq j\leq J-1}\|
x_{j+1}-x_j\|^q.
\end{equation}
\end{lem}
Denote $\mathcal{D}_0\subset \mathcal{D}_1\subset\ldots\subset
\mathcal{D}_r$ the algebras of intervals on $[0,1]$ obtained by
successive dyadic refinements. Define the $X$-valued function
$$\xi=\sum_{1\leq j\leq J-1}1_{[\frac{j}{J},\frac{j+1}{J}[\;}(x_{j+1}-x_j)$$
and consider expectations
$\xi_s=\mathbf{E}\left[\xi|\mathcal{D}_s\right]$ for
$s=1,\ldots,r$. Since $\xi_s$ form a martingale ranging in $X$, it
satisfies inequality (\ref{ineqlem1}). On the other hand
\begin{eqnarray*}
\| \xi_{s+1}-\xi_{s}\|_q^q & =& 2^{-r+s}2^{qs}\sum_{1\leq t\leq
2^{r-s}}^r2^{-qs}\|
2x_{t2^s}-x_{(t-1)2^s}-x_{(t+1)2^s}\|^q\\
                                         & \leq & 2^{-qs}
                                         \min_{2^s\leq j\leq
                                         J-2^s}\|
                                         2x_{j}-x_{(j-2^s}-x_{j+2^s}\|^q.
\end{eqnarray*}
So (\ref{ineqlem2}) follows from the fact that
\begin{eqnarray*}
\| \xi_{s}\|_q^q  \leq \| \xi_{s+1}-\xi_{s}\|_{\infty}^q
                                =  \sup_{j}\|
                                x_{j+1}-x_j\|^q.\;\blacksquare
\end{eqnarray*}

\begin{lem}
If $f_1,\ldots,f_J$, with $J=2^r$, is a finite system of functions
in $L^{\infty}_{X}(\Omega)$. Then
\begin{equation}\label{ineqlem3}
\sum_{s=1}^r2^{-qs}\min_{2^s\leq j\leq J-2^s} \|
2f_{j}-f_{j-2^s}-f_{j+2^s}\|^q\leq C\sup_{1\leq j\leq J-1}\|
f_{j+1}-f_j\|_{\infty}^q.
\end{equation}
\end{lem}
\bpr Replace $X$ by $L^q_X(\Omega)$, for which (\ref{ineqlem1})
remains valid, and use (\ref{ineqlem2}).\epr
\begin{lem}\label{lem4}
Let $f_1,\ldots,f_J$, with $J=2^r$, be a sequence of functions on
$\{1,-1\}^J$ where $f_j$ only depends on $\eps_1,\ldots, \eps_j$.
Then
\begin{eqnarray*}
\sum_{s=1}^r2^{-qs}\min_{2^s\leq j\leq J-2^s}\left(
\int_{\Omega_j\times \Omega_{2^s}\times \Omega_{2^s}}\|
f_{j+2^s}(\eps,\delta)-f_{j+2^s}(\eps,\delta')\|^q d\eps d\delta
d\delta'\right) \\ \leq 2^qC\sup_{1\leq j\leq J-1}\|
f_{j+1}-f_j\|_{\infty}^q.
\end{eqnarray*}
\end{lem}
\bpr For every $d\leq j\leq J-d$, using the triangle inequality,
we obtain
\begin{eqnarray*}
\| 2f_j-f_{j-d}-f_{j+d}\|^q_q & = & \int_{\Omega_j\times
\Omega_{d}}\|
 2f_j-f_{j-d}-f_{j+d}\|^q d\eps d\delta\\
                                            & \geq &
                                            2^{-q}\int_{\Omega_j\times
                                            \Omega_{d}\times
                                            \Omega_d}\|
                                            f_{j+2^s}(\eps,\delta)-f_{j+2^s}(\eps,\delta')\|^q
                                            d\eps
d\delta d\delta'.\\
\end{eqnarray*}
The lemma then follows from (\ref{ineqlem3}).\epr

\

Now, let us prove Theorem~\ref{embedtreethm}. Fix $J$ and consider
a $1$-Lipschitz map $F: T_J\to X$. Apply Lemma \ref{lem4} to the
functions $f_1,\ldots,f_J$ defined by
$$\forall \alpha\in
\Omega_j,\quad f_j(\alpha)=F(\alpha).$$

By definition of the compression, we have
\begin{equation}\label{ineq1}
\rho\left(d\left((\eps,\delta),(\eps,\delta)\right)\right)\leq \|
f_{j+2^s}(\eps,\delta)-f_{j+2^s}(\eps,\delta')\|
\end{equation}
where $\eps \in \Omega_j$ and $\delta,\delta'\in \Omega_{2^s}$.

But, on the other hand, with probability $1/2$, we have
$$d\left(\left(\eps,\delta),(\eps,\delta)\right)\right)=2.2^{s}.$$
So combining this with Lemma \ref{lem4}, (\ref{ineq1}) and with
the fact that $F$ is $1$-Lipschitz, we obtain
\begin{eqnarray*}
\sum_{s=1}^r2^{-qs}\rho(2^s)^q\leq 2^{q+1}C
\end{eqnarray*}
But since $\rho$ is decreasing, we have
$$2^{-qs}\rho(2^s)^q\geq 2^{-q-1}\int_{2^{s-1}+1}^{2^s}\frac{1}{t}\left(\frac{\rho(t)}{t}\right)^qdt.$$
So (\ref{thmineq}) follows.\epr

\section{Applications and further results}

\subsection{Hilbert compression, volume growth and random walks}

Let $G$ be a locally compact group generated by a symmetric
compact subset $S$ containing $1$. Let us denote $V(n)=\mu(S^n)$
and $S(n)=V(n+1)-V(n)=\mu(S^{n+1}\smallsetminus S^n)$. Extend $V$
as a piecewise linear function on $\R_+$ such that $V'(t)=S(n)$
for $t\in ]n,n+1[$.

\begin{prop}\label{croissanceinterm}
Let $G$ be a compactly generated locally compact group. For any
$2\leq p<\infty$,
$$J_{G,p}(t)\preceq \frac{t}{\log V(t)}.$$
\end{prop}
\bpr For every $n\in \N$, define $$k(n)=\sup\{k, \; V(n-k)\geq
V(n)/2\}$$ and
$$j(n)=\sup_{1\leq j\leq n}{k(j)}.$$
For every positive integer $l\leq n/j(n)$,
$$V(n) \geq 2^l V(n-lj(n))$$
Hence, as $V(0)=1,$
$$V(n)\geq 2^{n/(j(n)+1)}.$$
Thus, there is a constant $c>0$ such that
$$j(n)\geq \frac{cn}{\log V(n)}.$$
Let $q_n\leq n$ be such that $j(n)=k(q_n).$  Now define
$$\varphi_n=\sum_{k=1}^{q_n-1}1_{B(1,k)}.$$
Note that the subsets $SB(1,k)\vartriangle
B(1,k)=B(1,k+1)\smallsetminus B(1,k)$, for $k\in \N$, are
piecewise disjoint. Thus, an easy computation shows that
$$\|\tilde{\nabla} \varphi_n\|_p\leq V(q_n)^{1/p}.$$
On the other hand
$$\|\varphi_n\|_p\geq j(n)V(q_n-j(n))^{1/p}\geq \frac{cn}{\log V(n)}(V(q_n)/2)^{1/p}.$$
Since $J^b_{G,p}(n)\geq \|\varphi_n\|_p/\|\tilde{\nabla}
\varphi_n\|_p$, we conclude that $J^b_{G,p}(n)\succeq n/\log
V(n)$.\epr

\

Now, consider a symmetric probability measure $\nu$ on a finitely
generated group $G$, supported by a finite generating subset $S$.
Given an element $\varphi$ of $\ell^2(G)$, a simple calculation
shows that
$$\frac{1}{2}\int \int|\varphi(sx)-\varphi(x)|^2d\nu^{(2)}(s)d\mu(x)=\int (\varphi-\nu^{(2)}\ast \varphi)\varphi d\mu=\| \varphi\|_2^2-\|\nu\ast \varphi\|_2^2$$
where $\mu$ denotes the counting measure on $G$. Let us introduce
a (left) gradient on $G$ associated to $\nu$. Let $\varphi$ be a
function on $G$; define
$$|\tilde{\nabla}\varphi|^2_2(g)=\int |\varphi(sg)-\varphi(g)|^2d\nu^{(2)}(s).$$
This gradient satisfies
$$\||\tilde{\nabla}\varphi|_2\|_2^2=2(\| \varphi\|_2^2-\|\nu\ast \varphi\|_2^2).$$
We have
$$\mu(S)^{-1/2}|\tilde{\nabla}\varphi|_2\leq |\tilde{\nabla}\varphi|\leq |\tilde{\nabla}\varphi|_2.$$
\begin{prop}\label{marcheprop}
Assume that $\nu^{(n)}(1)\succeq e^{-Cn^b}$ for some $b<1$. Then
$$J^b_{G,2}(t)\succeq Ct^{1-b}.$$
\end{prop}
\bpr Let us prove that there exists a constant $C'<\infty$ such
that for every $n\in \N$, there exists $n\leq k\leq 2n$ such that
$$\frac{\parallel|\tilde{\nabla}\nu^{(2k)}|_2\parallel_2^2}{\parallel \nu^{(2k)}\parallel_2^2}\leq C'n^{b-1}.$$
Since $\nu^{(2k)}$ is supported in $S^{2k}\subset S^{4n}$, this
will prove the proposition. Let $C_n$ be such that for every
$n\leq q\leq 2n $
$$\frac{\parallel|\tilde{\nabla}\nu^{(2q)}|_2\parallel_2^2}{\parallel \nu^{(2q)}\parallel_2^2}\geq C_nn^{b-1}.$$
Since the function defined by $\psi(q)=\parallel
\nu^{(2q)}\parallel_2^2$ satisfies
$$\psi(q+1)-\psi(q)=-\frac{1}{2}\parallel|\tilde{\nabla}\nu^{(2q)}|_2\parallel_2^2,$$
we can extend $\psi$ as a piecewise linear function on $\R_+$ such
that
$$\psi'(t)=\frac{1}{2}\parallel|\tilde{\nabla}\nu^{(2q)}|_2\parallel_2^2$$
for every $t\in [q,q+1[$. Then, for every $n\leq t\leq 2n$ we have
$$-\frac{\psi'(t)}{\psi(t)}\geq C_nn^{b-1}$$
which integrates in
$$-\log\left(\frac{\psi(2n)}{\psi(n)}\right)\geq C_nn^{b}$$
Since $\psi(n)<1$, this implies
$$\psi(2n)\leq e^{-C_nn^{b}}.$$
But on the other hand, $$\psi(2n)\geq
\parallel\nu^{(4n)}\parallel_2^2\geq \nu^{(8n)}(1)\geq
e^{-8Cn^{b}}.$$ So $C_n\leq 8C.$\epr

\subsection{A direct construction to embed trees}

Here, we propose to show that the method used in
\cite{Bourgain,GK,BS} to embed trees in $L^p$-spaces can also be
exploited to obtain optimal estimates (i.e. a converse to
Theorem~\ref{embedtreethm}). Moreover, no hypothesis of local
finitude is required for this construction.

\begin{thm}\label{sphericalprop}
Let $T$ be a simplicial tree. For every increasing function
$f:\R_+\to\R_+$ satisfying, for $1\leq p<\infty$
\begin{equation}\tag{$C_p$}\quad
\int_1^{\infty}\left(\frac{f(t)}{t}\right)^p\frac{dt}{t}<\infty,
\end{equation}
there exists a uniform embedding $F$ of $T$ into $\ell^p(T)$ with
compression $\rho\succeq f.$
\end{thm}
\bpr Let us start with a lemma.

\begin{lem}
For every nonnegative sequence $(\xi_n)$ such that
$$\sum_n|\xi_{n+1}-\xi_n|^p<\infty,$$ there exists a Lipschitz map
$F:T\to \ell^p(T)$ whose compression $\rho$ satisfies
$$\forall n\in \N,\quad\rho(n)\geq \left(\sum_{j=0}^n\xi_j^p\right)^{1/p}.$$
\end{lem}
\bpr The following construction is a generalization of those
carried out in \cite{GK} and \cite{BS}. Fix a vertex $o$. For
every $y\in T$, denote $\delta_y$ the element of $\ell^p(T)$ that
takes value $1$ on $y$ and $0$ elsewhere. Let $x$ be a vertex of
$T$ and let $x_0=x,x_1\ldots, x_l=o$ be the minimal path joining
$x$ to $o$. Define
$$F(x)=\sum_{i=1}^l \xi_i\delta_{x_i}.$$
To prove that $F$ is Lipschitz, it suffices to prove that $\|
F(x)-F(y)\|_p$ is bounded for neighbor vertices in $T.$ So let $x$
and $y$ be neighbor vertices in $T$ such that
$d(o,y)=d(x,o)+1=l+1$. We have
$$\| F(y)-F(x)\|_p^p\leq \xi_0^p+\sum_{j=0}^l|\xi_{n+1}-\xi_n|^p.$$
On the other hand, let $x$ and $y$ be two vertices in $T$. Let $z$
be the last common vertex of the two geodesic paths joining $o$ to
$x$ and $y$. We have $$d(x,y)=d(x,z)+d(z,y)$$ and
\begin{eqnarray*}
\| F(x)-F(y)\|_p^p & = & \| F(x)-F(z)\|_p^p+\| F(z)-F(y)\|_p^p\\
& \geq & \max\{\| F(x)-F(z)\|_p^p\; ,\;\| F(z)-F(y)\|_p^p\}.
\end{eqnarray*}
Let $k=d(z,x)$; we have
\begin{eqnarray*}
\| F(x)-F(z)\|_p^p & \geq & \sum_{j=0}^k \xi_j^{p},
\end{eqnarray*}
which proves the lemma.\epr

\

Now, let us prove the proposition. Define $(\xi_j)$ by
$$\xi_0=\xi_1=0;$$
$$\forall j\geq 1,\quad \xi_{j+1}-\xi_j=\frac{1}{j^p}\frac{f(j)}{j}$$
and consider the associated Lipschitz map $F$ from $T$ to
$\ell^p(T)$. Clearly, we have
$$\sum |\xi_{n+1}-\xi_n|^p<\infty$$
and
\begin{eqnarray*}
\sum_{j=0}^n\xi_j^p & \geq &
\sum_{j=[n/2]}^n\left(\sum_{k=0}^{j-1}|\xi_{k+1}-\xi_k|\right)^p\\
& \geq & n/2\left(\sum_{k=0}^{[n/2]-1}|\xi_{k+1}-\xi_k|\right)^p\\
& \geq & cf([n/2])
\end{eqnarray*}
using the fact that $f$ is nondecreasing. So the proposition now
follows from the lemma.\epr

\subsection{Cocycles with lacunar compression}

\begin{prop}\label{Cpbehavior}
For any increasing sublinear function $h:\R_+\rightarrow \R_+$ and
every $2\leq p<\infty$, there exists a function $f$ satisfying
$(C_p)$, a constant $c>0$ and a increasing sequence of integers
$(n_i)$ such that
$$\forall i\in \N, f(n_i)\geq ch(n_i).$$
\end{prop}
\bpr Choose a sequence $(n_i)$ such that
$$\sum_{i\in \N}\left(\frac{h(n_i)}{n_i}\right)^p<\infty$$
Define
$$\forall i\in \N, n_i\leq t<n_{i+1},\quad f(t)=h(n_i)$$
We have
\begin{eqnarray*}
\int_{1}^{\infty} \frac{1}{t}\left(\frac{f(t)}{t}\right)^pdt &
\leq & \sum_i
(h(n_i))^p\int_{n_i}^{n_{i+1}}\frac{dt}{t^{p+1}}\\
& \leq & (p+1)\sum_i \left(\frac{h(n_i)}{n_i}\right)^p\\
& < & \infty
\end{eqnarray*}
So we are done. \epr

\subsection{The case of $\Z\wr\Z$}

The proof of Theorem~\ref{thmZwrZ} follows from
Proposition~\ref{Wp/compression} and from the following
observation.

\begin{prop}\label{ZwrZProp}
For all $1\leq p<\infty$, the maximal $\ell^p$-compression
function of the group $G=\Z\wr\Z$ satisfies
$$M\rho_{G,p}(t)\succeq t^{p/(2p-1)}.$$
\end{prop}
\bpr Denote by $\theta$ the projection $\Z\wr\Z\to C_2\wr\Z$. Fix
two word lengths on $\Z\wr\Z$ and $C_2\wr\Z$, which for
simplicity, we will both denote by $|g|$.

Consider the unique cocycle $b:\Z\wr\Z\to \ell^p(\Z)$ which
extends the natural injective morphism $\Z^{(\Z)}\to\ell^p(\Z)$.
For any $g=(k,u)\in \Z\wr\Z=\Z\ltimes\Z^{(\Z)}$, we therefore have
$\|b(g)\|=\|u\|_p.$ Taking the $\ell^p$-direct sum of this cocycle
with every cocycle of $\Z\wr\Z$ factorizing through $\theta$, and
since $M\rho_{C_2\wr\Z,p}(t)\approx t$, we obtain
\begin{equation}\label{MaxCompEq}
M\rho_{\Z\wr\Z,p}(t)\succeq \inf_{g\in \Z\wr\Z,\;|g|\geq
t}\max\{|p(g)|,\|b(g)\|\}.
\end{equation}

Up to multiplicative constants, (see \cite[Theorem 1.2]{Par}), the
word length of an element $g=(k,u)\in \Z\wr\Z$ is given by
$$L(\gamma)+\sum_{h\in \Z}|u(h)|=L(\gamma)+\|u\|_1,$$ where
$L(\gamma)$ is the length of a shortest path $\gamma$ from $0$ to
$k$ passing through every element of the support of $u$.
Similarly, $|p(g)|\approx L(\gamma)+|\Supp(u)|.$ Hence by
(\ref{MaxCompEq}), we can assume that $L(\gamma)\leq |g|/2$, so
that $\|u\|_1\geq |g|/2$. By Hölder's inequality, we have
$\|u\|_1\leq \|u\|_p|\Supp(u)|^{1-1/p}$, which is less than a
constant times $\|b(g)\||p(g)|^{1-1/p}.$ Therefore
$$M\rho_{\Z\wr\Z,p}(t)\succeq \inf_{g\in \Z\wr\Z,\;|g|\geq
t}\max\left\{|p(g)|,|g|/|p(g)|^{1-1/p}\right\},$$ which
immediately implies the proposition. \epr

\bigskip
\footnotesize

\noindent Romain Tessera\\
\noindent Equipe Analyse, Géométrie et Modélisation\\
Université de Cergy-Pontoise, Site de Saint-Martin\\
2, rue Adolphe Chauvin F 95302 Cergy-Pontoise Cedex, France\\
E-mail: \url{tessera@clipper.ens.fr}\\

\end{document}